\title{On two biased graph processes}
\author{{Gideon Amir \thanks{Weizmann Institute, Rehovot, 76100,
Israel. Email: gideon.amir@weizmann.ac.il}} \quad {Eyal Lubetzky
\thanks{ School of Computer Science, Raymond and Beverly
Sackler Faculty of Exact Sciences, Tel Aviv University, Tel Aviv,
69978, Israel. Email: lubetzky@tau.ac.il. Research partially
supported by a Charles Clore Foundation Fellowship.}} }

\documentclass [11pt]{article}
\usepackage[]{graphicx,epsf,epsfig,amsmath,amssymb,amsfonts,latexsym,amsthm}

\newtheorem{theorem}{Theorem}[section]
\newtheorem{lemma}[theorem]{Lemma}
\newtheorem{claim}[theorem]{Claim}
\newtheorem*{definition}{Definition}

\renewcommand{\epsilon}{\varepsilon}

\newcommand{\Gorig}[1][1]{{\mathcal{G}_#1}}

\newcommand{\Gor}[1][K]{{\mathcal{G}^{\vee}_#1}}
\newcommand{\Gorap}[1][K]{{\widetilde{\mathcal{G}^{\vee}_#1}}}
\newcommand{\Gand}[1][K]{{\mathcal{G}^{\wedge}_#1}}
\newcommand{\Gandap}[1][K]{{\widetilde{\mathcal{G}^{\wedge}_#1}}}
\newcommand{\tgand}{t^{\wedge}_g}
\newcommand{\tgor}{t^{\vee}_g}
\newcommand{\tcand}{t^{\wedge}_c}
\newcommand{\tcor}{t^{\vee}_c}

\newtheoremstyle{upright}%
        {8pt plus2pt minus4pt}%
        {8pt plus2pt minus4pt}%
        {\upshape}%
        {}%
        {\bfseries\scshape}%
        {:}%
        {1em}%
        {}%
\theoremstyle{upright}

\newcommand{\ds}{\displaystyle}
\newcommand{\ignore}[1]{}

\begin{document}
\maketitle

\begin{abstract}

In \cite{AGLS}, the authors consider the generalization $\Gor$ of
the Erd\H{o}s-R\'enyi random graph process $\Gorig$, where instead
of adding new edges uniformly, $\Gor$ gives a weight of size $1$ to
missing edges between pairs of isolated vertices, and a weight of
size $K\in[0,\infty)$ otherwise. This can correspond to the linking
of settlements or the spreading of an epidemic. The authors
investigate $\tgor(K)$, the critical time for the appearance of a
giant component as a function of $K$, and prove that
$t_g^\vee=\left(1 + o(1)\right)\frac{4}{\sqrt{3K}}$, using a proper
timescale.

In this work, we show that a natural variation of the model $\Gor$
has interesting properties. Define the process $\Gand$, where a
weight of size $K$ is assigned to edges between pairs of
non-isolated vertices, and a weight of size $1$ otherwise. We prove
that the asymptotical behavior of the giant component threshold is
essentially the same for $\Gand$, and namely $\tgand / \tgor$ tends
to $\frac{64\sqrt{6}}{\pi(24+\pi^2)}\approx 1.47$ as $K\to\infty$.
However, the corresponding thresholds for connectivity satisfy
$\tcand / \tcor=\max\{\frac{1}{2},K\}$ for every $K>0$. Following
the methods of \cite{AGLS}, $\tgand$ is characterized as the
singularity point to a system of differential equations, and
computer simulations of both models agree with the analytical
results as well as with the asymptotic analysis.
%
In the process, we answer the following
 question: when
does a giant component emerge in a graph process where edges are
chosen uniformly out of all edges incident to isolated vertices,
while such exist, and otherwise uniformly? This corresponds to the
value of $\tgand(0)$, which we show to be
$\frac{3}{2}+\frac{4}{3\mathrm{e}^2-1}$.
\end{abstract}

\section{Introduction}\label{sec::intro}
\subsection{The Erd\H{o}s-R\'enyi graph process and biased
processes} The random graph process, $\Gorig[1](n)$, is a sequence
of $\binom{n}{2}+1$ graphs on $n$ vertices,
$\Gorig[1]^0,\ldots,\Gorig[1]^{\binom{n}{2}}$, where $\Gorig[1]^0$
is the edgeless graph on $n$ vertices, and $\Gorig[1]^T$ is obtained
by adding an edge to $\Gorig[1]^{T-1}$, chosen uniformly over all
missing edges. This model was introduced by Erd\H{o}s and R\'enyi in
\cite{ErdosRenyi}, where it is shown that for every constant
$C<\frac{1}{2}$, the largest component of $\Gorig[1]^{C n}$ is
typically of size $O(\log n)$, and yet for every constant
$C>\frac{1}{2}$ there is typically a single component of linear size
in $\Gorig[1]^{Cn}$ and all other components are of size $O(\log
n)$. This single component and its evolution are referred to as the
``giant component'', and ``the double jump phenomenon''
respectively. Using a timescale of $n/2$ edges, the threshold for
the appearance of the giant component is thus $t_g = 1$. Another
classical result of Erd\H{o}s and R\'enyi determines that
$\Gorig[1]^{C n\log n}$ is typically connected for $C >
\frac{1}{2}$, and typically disconnected for $C < \frac{1}{2}$.
Using a timescale of $\frac{n}{2}\log n$, the connectivity threshold
is thus $t_c =1$. For further information on the evolution of the
giant component in the random graph process, as well as on the
threshold for connectivity, see, e.g., \cite{RandomGraphs}.

There has been extensive study on the thresholds for the appearance
of a giant component and for connectivity in different variations of
the random graph process. For instance, in a model suggested by
Achlioptas, two random edges are chosen uniformly out of the missing
edges at each step, out of which some algorithm $\mathcal{A}$
selects one to be added to the graph. For results on upper and lower
bounds on the emerging of the giant components for various
algorithms in this model, see
\cite{BohmanFrieze},\cite{BohmanKim},\cite{BohmanKravitz},\cite{SpencerWormald}.

In \cite{KangKohLuczakSangwook}, the authors consider a random graph
process on multi-graphs, where at each step an edge is added between
a random vertex of minimal degree and a random uniformly chosen
vertex. The authors analyze the number of vertices of degrees
$0,1,2$ along the process using the differential equation method for
graph processes of Wormald \cite{WormaldDiffEq}, and show that the
mentioned graph process becomes connected typically when the minimal
degree becomes $3$.

The following generalization of the original graph process $\Gorig$,
which we denote by $\Gor(n)$, was studied in \cite{AGLS}: at each
step, $\Gor$ gives a weight of size $1$ to missing edges between
pairs of isolated vertices, and a weight of size $K\in[0,\infty)$
otherwise (when no isolated vertices are left, the distribution on
the missing edges becomes the uniform distribution). This can
correspond to the linking process of $n$ initially isolated
settlements, or the spreading of an epidemic, where the probability
of a new link is affected by whether or not one of its endpoints
already has other links. The threshold for the appearance of a giant
component in $\Gor$, $\tgor$, becomes a continuous function of the
parameter $K$, and the authors of \cite{AGLS} use the differential
equation method to express $\tgor$ as a singularity point to a
system of coupled non-linear ordinary differential equations (ODEs).
By applying methods of asymptotic analysis of ODEs, it is proved
that $\tgor(K)=\left(1+o(1)\right)\frac{4}{\sqrt{3K}}$, where the
$o(1)$-term tends to $0$ as $K\to\infty$.

In this work, we show that a natural variation on the process $\Gor$
has interesting properties. Consider the process $\Gand$, where
instead of placing the weight $K$ when one of the endpoints is
non-isolated (as in $\Gor$), it is placed when {\em both} endpoints
are non-isolated. In other words, $\Gor$ gives a weight of size
$K\in[0,\infty)$ to missing edges between pairs of non-isolated
vertices, and a weight of size $1$ otherwise. Notice that for $K=1$,
both processes are equivalent to the original graph process.
Furthermore, for any $K$, both processes $\Gor$ and $\Gand$ apply
the rule of the original graph process $\Gorig$ once no isolated
vertices are left, hence it is interesting to compare the two up to
roughly that point.

Applying the methods of \cite{AGLS} on $\Gand$, we express its
threshold for the appearance of a giant component, $\tgand(K)$, as a
singularity point to a system of coupled non-linear ODEs. For the
special case $K=0$, we obtain a setting where edges are added
uniformly at random out of all edges incident to an isolated vertex,
until no such vertex is left (from that point on, edges are added
uniformly at random). In this special case we get
$\tgand(0)=\frac{3}{2}+\frac{4}{3\mathrm{e}^2-1}$.

Using classical methods of asymptotic analysis of ODEs (see, e.g.,
\cite{Bender}) we readily calculate the asymptotic behavior of
$\tgand(K)$ for $K \gg 1$, and obtain that
$\tgand(K)=\left(1+o(1)\right)\frac{\pi}{2\sqrt{2}}\left(1+\frac{\pi^2}{24}\right)\frac{1}{\sqrt{K}}$.
 It follows that $\tgand / \tgor \approx 1.47$ for $K
\gg 1$, and we note that obtaining this result via combinatorial
arguments seems challenging. Numerical approximations of the ODEs
validate this asymptotic analysis.

While the behavior of the threshold for the appearance of a giant
component is similar for $K \gg 1$, combinatorial arguments yield
that for all $K>0$, $\tcor(K) = 1$, whereas
$\tcand(K)=\max\{\frac{1}{2},K\}$, hence $\tcand / \tcor = K$ for
every $K\geq \frac{1}{2}$.

It is possible to implement both processes $\Gand$ and $\Gor$
efficiently by choosing an appropriate data-structure for holding
the isolated and non-isolated vertices. Our implementation requires
$O(n)$ memory and runs in time $O(n\log n)$, and its results
validate the above analytical results.

\subsection{Notations and main results}
A property of graphs is a collection of graphs closed under
isomorphism. A property is said to be monotone (increasing) if it is
closed under the addition of edges. Throughout the paper, we say
that a property of graphs on $n$ vertices occurs {\em with high
probability}, or {\em almost surely}, or that {\em almost every}
process $\mathcal{G}$ satisfies this property, if the probability
for the corresponding event tends to $1$ as $n\to \infty$. Note
that, when proving that certain statements hold with high
probability, one may clearly condition on events which hold with
high probability.

On several occasions, we examine the processes
 $\Gor$ or $\Gand$ starting from some arbitrary graph $H$ (instead of
the edgeless graph). We denote these processes by $\Gor|_H$ and
$\Gand|_H$.

Given a process $\mathcal{G}$, we let $\mathcal{G}^i$ denote
$\mathcal{G}$ after $i$ edges. It will be convenient to use two
timescales when referring to $\mathcal{G}$, which we denote by:
\begin{eqnarray}\mathcal{G}(t) &:=& \mathcal{G}^{tn/2}~,\nonumber\\
\mathcal{G}[t] &:=& \mathcal{G}^{t\frac{n}{2}\log
n}~.\nonumber\end{eqnarray} We say that $t$ is a threshold for the
appearance of a giant component in $\mathcal{G}$ if for every
$\epsilon>0$, with high probability $\mathcal{G}(t-\epsilon)$ does
not contain a giant component and yet $\mathcal{G}(t+\epsilon)$ does
contain one. Similarly, we say that $t$ is a threshold for
connectivity if for every $\epsilon>0$, with high probability
$\mathcal{G}[t-\epsilon]$ is disconnected whereas
$\mathcal{G}[t+\epsilon]$ is connected. According to these
definitions, both thresholds equal $1$ for the original random graph
process $\Gorig$.

Theorem \ref{thm-conn} determines the threshold for connectivity in
both processes, $\tcor$ and $\tcand$, and shows that $\tcor = 1$ for
any $K>0$ whereas $\tcand= K$ for any $K > \frac{1}{2}$, hence their
ratio is $K$ for any $K>\frac{1}{2}$. Furthermore, as we later
state, it follows from \cite{AGLS} that a ratio of $\lceil
\max\{\frac{1}{K},K\}\rceil$ is the maximal possible between the
threshold of $\Gand$ and the threshold of $\Gorig$ for any monotone
property, and therefore, $\tcor$ achieves this maximum.
\begin{theorem}\label{thm-conn}
For every $K > 0$, $\tcor(K)=1$ and yet
$\tcand(K)=\max\{\frac{1}{2},K\}$. In the special case $K=0$ we have
$\tcor(0)=\tcand(0) = \frac{1}{2}$.
\end{theorem}

For a given graph $G=(V,E)$ on $|V|=n$ vertices, we let
$\mathcal{C}=\mathcal{C}(G)$ denote the set of connected components
of $G$, and for $i\in\mathbb{N}$, we define
$\mathcal{C}_i=\mathcal{C}_i(G)$ to be the set of components of size
$i$: $ \mathcal{C}_i=\{ C\in\mathcal{C}(G):|C|=i\}$. Whenever it is
clear as to which graph process $\mathcal{G}$ we are referring, we
use the abbreviation $\mathcal{C}_i^t$ to denote
$\mathcal{C}_i(\mathcal{G}^t)$. The fractions of vertices which
belong to components of size $1$ and $2$ are defined as:
$$ I(G) = \frac{|\mathcal{C}_1|}{n}~,~I_2(G)=\frac{2|\mathcal{C}_2|}{n}~,$$
and the {\em susceptibility} of $G$, $S(G)$, is defined to be the
average size of a connected component, averaged over all vertices:
$$S(G) = \frac{1}{n}\sum_{v \in V} |C(v)|=\frac{1}{n}\sum_{C \in \mathcal{C}(G)} |C|^2 ,$$ where $C(v)$ denotes the
connected component of $v$. The relation between $S(G)$ and the
existence of a giant component in $G$ is immediate: if $G$ contains
a component of size $\alpha n$ for some $\alpha > 0$, then $S(G)
\geq \alpha^2 n$, and if $S(G) \geq \alpha n$ then clearly there
exists a component of size $\alpha n$. Therefore, $G$ has a giant
component iff $S(G)=\Theta(n)$.

The methods used in \cite{AGLS} to analyze $\Gor$ and describe
$\tgor$ as a singularity point to a system of ODEs can in fact be
applied to a wider class of graph processes. Namely, these methods
can be applied to every process where the weight function $W$ on the
missing edges satisfies $\max W / \min W \leq K$ for some constant
$K>0$. Instead of repeating the complete set of arguments of
\cite{AGLS} in order to prove analogous results on $\Gand$, we
summarize the arguments briefly, and proceed to prove the necessary
conditions required for them to work.

Following the ideas of \cite{AGLS} and previously of
\cite{SpencerWormald}, define the following system of coupled ODEs:
\begin{eqnarray}
&\left\{\begin{array}{lcl}
y' &=& \displaystyle{\frac{-y}{1+(K-1)(1-y)^2}}\\
y(0)&=&1 \\
\end{array} \right. ~,\label{y-eq}\\
&\left\{\begin{array}{lcl}
z' &=& \displaystyle{\frac{z^2+(K-1)(z-y)^2}{1+(K-1)(1-y)^2}}\\
z(0)&=&1 \\
\end{array} \right. ~.\label{z-eq}
\end{eqnarray}
We further define: \begin{equation}
\label{w-eq}\left\{\begin{array}{lcl}
w' &=& \displaystyle{\frac{y^2-2wy -2Kw(1-y)}{1+(K-1)(1-y)^2}}\\
w(0)&=&0 \\
\end{array} \right.~,
\end{equation} where $y$ is the solution to \eqref{y-eq}.
Theorem \ref{thm-tgand} states that $I(G)$, $S(G)$ and $I_2(G)$
along the process $\Gand$ are approximated by the solutions to the
ODEs above, and that $\tgand$ is equal to the singularity point of
the solution to \eqref{z-eq} for any $K>0$:
\begin{theorem}\label{thm-tgand}
Let $y(t)$, $z(t)$ and $w(t)$ denote the solutions for
\eqref{y-eq},\eqref{z-eq},\eqref{w-eq}, and let $x_c$ denote the
singularity point of $z(t)$ if such exists, and $\infty$ otherwise.
For $0 < \delta < 1$, let $\tau_\delta > 0$ be the minimal point
satisfying $y(\tau_\delta)\leq \delta$. The following statements
hold almost surely:
\begin{enumerate}
\item For every $\delta > 0$,
$\left|I\left(\Gand(t)\right)-y(t)\right|=o(1)$ and
$\left|I_2\left(\Gand(t)\right)-w(t)\right|=o(1)$ for all $0 \leq t
\leq \tau_\delta$.
\item For every $\epsilon,\delta >0$, $\left|S\left(\Gand(t)\right)-z(t)\right|=o(1)$
for all $0\leq t \leq \min\{\tau_\delta, x_c -\epsilon\}$.
\item For all $K>0$, $\tgand = x_c$.
\end{enumerate}
\end{theorem}
The value of $\tgand$ in the special case $K=0$ and the asymptotic
behavior of $\tgand$ are stated in the following theorem:
\begin{theorem}\label{thm-K-0-K-gg-1}
The threshold for the appearance of a giant component, $\tgand(K)$,
is a continuous function of $K$ on $(0,\infty)$, and satisfies:
$$\left\{\begin{array}{lcl}
\tgand(0)&=&\frac{3}{2}+\frac{4}{3\mathrm{e}^2-1}\\
\tgand(K) &=&\left(1 + o(1)\right)
\frac{\pi}{2\sqrt{2}}\left(1+\frac{\pi^2}{24}\right)\frac{1}{\sqrt{K}}\\
\end{array} \right.~,$$
where the $o(1)$-term tends to $0$ as $K\to\infty$.
\end{theorem}

The rest of this paper is organized as follows: in Section
\ref{sec::conn}, we prove Theorem \ref{thm-conn} by examining second
moments of processes which are easier to analyze, and provide
results of computer simulations of $\tcor$ and $\tcand$. In Section
\ref{sec::giant}, we prove Theorems \ref{thm-tgand} and
\ref{thm-K-0-K-gg-1} using analysis of differential equations, and
provide results of computer simulations of $\tgand$.

\section{Thresholds for connectivity}\label{sec::conn}
The proof of Theorem \ref{thm-conn} uses the following result of
\cite{AGLS}:

\begin{definition} Let $M \in \mathbb{N}$.
An \texttt{$M$-bounded weighted graph process} on $n$ vertices,
$\mathcal{H}=\mathcal{H}(n)$, is an infinite sequence of graphs on
$n$ vertices, $(\mathcal{H}^0,\mathcal{H}^1,\ldots)$, where
$\mathcal{H}^0$ is some fixed initial graph, and $\mathcal{H}^t$ is
generated from $\mathcal{H}^{t-1}$ by adding one edge at random,
according to a distribution of the following type: the probability
of adding the edge $e$ to $\mathcal{H}^{t-1}$ is proportional to
some weight function $W_t(e)$, satisfying:
$$\max_{e\notin \mathcal{H}^{t-1}} W_t(e) \leq M \min_{e\notin H^{t-1}}
W_t(e)~.
$$If for some $\nu\geq 0$ $\mathcal{H}^\nu = K_n$, we define
$\mathcal{H}^t = \mathcal{H}^\nu = K_n$ for every $t>\nu$.
\end{definition}

\begin{theorem}[\cite{AGLS}]\label{thm-M-bounded}
Let $\mathcal{H}$ denote an $M$-bounded weighted graph process on
$n$ vertices, and let $\mathcal{A}$ denote a monotone increasing
property of graphs on $n$ vertices. The following statements hold
for any $t \in \mathbb{N}$:
\begin{eqnarray}
\label{dom-M-by-1} \Pr[\mathcal{H}^t \in \mathcal{A}]
&\leq& \Pr[\Gorig^{M t}|_{\mathcal{H}^0} \in \mathcal{A} ] ~,\\
\label{dom-1-by-M} \Pr[\Gorig^{t}|_{\mathcal{H}^0} \in
\mathcal{A}]&\leq& \Pr[\mathcal{H}^{M t}\in \mathcal{A}]~.
\end{eqnarray}
\end{theorem}

Notice that for every $K>0$, the processes $\Gor$ and $\Gand$ are
both $M$-bounded, where $M=\lceil\max\{\frac{1}{K},K\}$. As it is
well-known that $I(t)$, the fraction of isolated vertices is
$(1+o(1))\exp(-C)$ at time $Cn/2$, Theorem \ref{thm-M-bounded}
implies that the for any fixed $C$ and $K>0$, $\Gor(C)$ and
$\Gand(C)$ almost surely contain $\Theta(n)$ isolated vertices. On
the other hand, it is not difficult to see that $\Gor[0](1)$ has at
most $1$ isolated vertex, and we will later show that
$\Gand[0](\frac{3}{2}+o(1))$ almost surely has no isolated vertices.
For this reason, we treat the cases $K>0$ and $K=0$ separately.
Lemmas \ref{lem-conn-or-model} and \ref{lem-conn-and-model} below
establish the precise behavior of these parameters when $K>0$:

\begin{lemma}\label{lem-conn-or-model}
Let $K>0$. For every $0 < \epsilon < 1$, $\Gor[K][1+\epsilon]$ is
almost surely connected, whereas $\Gor[K][1-\epsilon]$ is almost
surely disconnected. Altogether, $\tcor=1$.
\end{lemma}
\begin{lemma}\label{lem-conn-and-model}
Let $K>0$. For every $0 < \epsilon < 1$,
$\Gand[K][(1+\epsilon)\max\{\frac{1}{2},K\}]$ is almost surely
connected, whereas $\Gand[K][K-\epsilon]$ and
$\Gand[K][\frac{1}{2}-\epsilon]$ are almost surely disconnected.
Altogether, $\tcand=\max\{\frac{1}{2},K\}$.
\end{lemma}
As we stated in the introduction, combining the fact that $t_c=1$
for $\Gorig$ with Theorem \ref{thm-M-bounded} yields that $\tcand(K)
\leq \lceil \max\{\frac{1}{K},K\}\rceil$, and indeed Lemma
\ref{lem-conn-and-model} shows that for $K\geq 1$, $K \in
\mathbb{N}$, this maximum is achieved.

The remaining case $K=0$ is settled by the next lemma:
\begin{lemma}\label{lem-conn-K-0}
The threshold for connectivity for $K=0$ satisfies
$\tcor(0)=\tcand(0)=\frac{1}{2}$.
\end{lemma}

\subsection{Proof of Lemma \ref{lem-conn-or-model}}
Let $0 < \epsilon < 1$; we first show that $\tcor \leq 1+\epsilon$.
It is easy and well known that for every $\epsilon'>0$, there exists
some $c=c(\epsilon')$ such that, with high probability, $\Gorig(c)$
has a giant component of size at least $(1-\epsilon')n$: assume that
indeed this holds. Take $0 < \epsilon' <
\ds{\min\{\frac{\epsilon/4}{1+\epsilon},K\}}$ and $c=c(\epsilon')$
as above, and notice that:
\begin{equation} \frac{1-\epsilon'}{1+\epsilon'} \geq
\frac{1+\frac{\epsilon}{2}}{1+\epsilon}~.
\end{equation}
By Theorem \ref{thm-M-bounded}, the giant component of
$\Gor^{c'n/2}$ typically contains at least $(1-\epsilon')n$
vertices, where $c' = \lceil \max\{K,\frac{1}{K}\} \rceil c$. Let
$H^i$ denote the largest component of $\Gor^i$, and let $A_v^i$, $v
\in V$, denote the event: $(v \notin H^i)$. The following holds for
every $i \geq c' n/2$:
\begin{equation}\label{eq-or-av-prob}
\Pr[\neg A_v^{i+1} | A_v^i] = \frac{ K|H^i| }{K (\binom{n}{2}-i) -
(K-1)\binom{|\mathcal{C}_1^i|}{2}} \geq \frac{K
(1-\epsilon')n}{K\binom{n}{2}+ \binom{|\mathcal{C}_1^i|}{2}} \geq
\frac{1-\epsilon'}{(1+\epsilon')\frac{n}{2}} ~,\end{equation} where
the last inequality is by the fact that $\frac{|\mathcal{C}_1^i|}{n}
\leq \epsilon' < K$. Setting $T_0 = c' n/2$, this gives the
following upper bound on $A_v^T$ for $v \notin H^{T_0}$:
$$ \Pr[A_v^T | \neg A_v^{T_0}] = \prod_{j=T_0}^T \Pr[A_v^j | A_v^{j-1}]
\leq \prod_{j=T_0}^T
\left(1-\frac{1-\epsilon'}{(1+\epsilon')\frac{n}{2}} \right) \leq
$$
\begin{equation}\label{eq-or-av-bound} \leq \exp\left(-(T-T_0)
\frac{1-\epsilon'}{(1+\epsilon')\frac{n}{2}}\right)~.
 \end{equation}
Thus, for $T = (1+\epsilon)\frac{n}{2}\log n$ we get:
$$ \Pr[A_v^T] \leq \exp\left( -\frac{1-\epsilon'}{1+\epsilon'}\left((1+\epsilon)\log n - c'\right)\right) \leq
\exp\left( -(1+\frac{\epsilon}{2})\log n + O(1)\right) =
o(n^{-1})~,$$ hence a union bound on the vertices of $V \setminus
H^{T_0}$ implies that, with high probability, $\Gor[K][1+\epsilon]$
is connected.

The lower bound on $\tcor$ is slightly more delicate. A simple way
to show the lower bound in $\Gorig$ is to apply a second moment
argument on the number of isolated vertices of $\Gorig(1)$. However,
in the case of $\Gor$, using uniform upper and lower bounds on
$\mathcal{C}_1^i$ (such as $\epsilon' n$ and $o(n)$) does not yield
useful bounds on the above second moment. The following claim
resolves this difficulty:
\begin{claim}\label{clm-iso-exponent-or} Let $0 < \delta < \frac{1}{2}$, and let
$\{H_0\}_n$ denote a family of arbitrary graphs on $n$ vertices with
$(1+o(1))n^{1-\delta}$ isolated vertices. If $\mathcal{H}$ is a
biased $\Gor$ process on $n$ vertices which begins with $H_0$, that
is: $\mathcal{H} \sim \Gor|_{H_0}$, then
 $\mathcal{H}[1-2\delta]$ almost surely contains isolated vertices.
\end{claim}
Indeed, to obtain the lower bound on $\tcor$ from the above claim,
fix $\delta = \epsilon/2$, and let $\tau$ denote the minimal time at
which $\mathcal{C}_1^\tau \leq n^{1-\delta}$. Taking
$H_0=\Gor^\tau$, Claim \ref{clm-iso-exponent-or} implies that
$\Gor|_{H_0}$ almost surely has isolated vertices at time $\tau +
(1-2\delta)\frac{n}{2}\log n$, and in particular, we obtain that
$\tcor \geq 1-\epsilon$, as required. It remains to prove Claim
\ref{clm-iso-exponent-and}:
\begin{proof}[Proof of claim] Let $H_0$ and $\delta$ be as above, and let $\mathcal{H} \sim
\Gor|_{H_0}$. Let $B_u^i$ denote the event that the vertex $u$ is
isolated at time $i$, where $i \leq n \log n$. Clearly:
$$
\Pr[\neg B_u^{i+1} | B_u^i]  = \frac{ K (n-|\mathcal{C}_1^i|) +
(|\mathcal{C}_1^i|-1) }{K (\binom{n}{2}-i) -
(K-1)\binom{|\mathcal{C}_1^i|}{2}} ~.$$ By the assumption on $H_0$,
$|\mathcal{C}_1^i| \leq |\mathcal{C}_1^0| = (1+o(1))n^{1-\delta}$,
thus:
\begin{equation}\label{eq-or-pr-b-u}\Pr[\neg B_u^{i+1} | B_u^i] =
\frac{2+O(n^{-\delta})}{n\left(1+O(\frac{\log n}{n}) +
O(n^{-2\delta})\right)}= \frac{2+O(n^{-\delta})}{n}~.\end{equation}
Similarly, we can define $B_{u,v}^i=B_u^i \wedge B_v^i$ for $i\leq
n\log n$ and get:
\begin{equation}\label{eq-or-pr-b-u-v}
\Pr[\neg B_{u,v}^{i+1} | B_{u,v}^i] = \frac{ 2K
(n-|\mathcal{C}_1^i|) + 2(|\mathcal{C}_1^i|-2)+1 }{K
(\binom{n}{2}-i) - (K-1)\binom{|\mathcal{C}_1^i|}{2}} =
\frac{4+O(n^{-\delta})}{n}~.\end{equation} Let
$T=(1-2\delta)\frac{n}{2}\log n$, and define $Y =
\sum_{u}\mathbf{1}_{B_u^T} = |\mathcal{C}_1^T|$ to be the number of
isolated vertices of $\mathcal{H}^T$. A straightforward second
moment consideration implies that $Y
> 0$ almost surely. To see this, first note that
\eqref{eq-or-pr-b-u} along with the well known bound
$1-x\geq\mathrm{e}^{-x/(1-x)}$ for $0\leq x < 1$ yield:
$$
\mathbb{E}Y=|\mathcal{C}_1^0| \prod_{i=1}^T \Pr[B_u^i | B_u^{i-1}] =
(1+o(1))n^{1-\delta} \left(1-\frac{2+O(n^{-\delta})}{n}\right)^T
\geq $$
\begin{equation}\label{eq-or-Y-exp}\geq (1+o(1))n^{1-\delta}
\exp\left(-\frac{\left(2+o(1)\right)/n}{1-\frac{2+o(1)}{n}}
(1-2\delta)\frac{n}{2}\log n\right) = n^{1-\delta}
n^{-1+2\delta+o(1)} = n^{\delta+o(1)}~.
\end{equation}
By \eqref{eq-or-pr-b-u} and \eqref{eq-or-pr-b-u-v}, there exist
$p=p(n)$ and $q=q(n)$ such that $\Pr[\neg B_{u,v}^{i+1} | B_{u,v}^i]
\geq p = \frac{4+O(n^{-\delta})}{n}$ and $\Pr[\neg B_u^{i+1} |
B_u^i] \leq q = \frac{2+O(n^{-\delta})}{n}$. The following holds:
$$ \mathrm{Cov}(\mathbf{1}_{B_u^T},\mathbf{1}_{B_v^T}) = \Pr[B_{u,v}^T] - \Pr[B_{u}^T]^2
\leq (1-p)^T - (1-q)^{2T} =
$$
$$ =
\left(1-p-(1-q)^2\right) \sum_{i=0}^{T-1}(1-p)^i(1-q)^{2(T-1-i)} =
\frac{O(n^{-\delta})}{n} T \left(1-\frac{4+o(1))}{n}\right)^{T-1}
\leq
 $$
$$ \leq (1-2\delta) \log n \cdot O(n^{-\delta}) \exp\left(-(1+o(1))\frac{4}{n}(1-2\delta)\frac{n}{2}\log n
\right) = n^{-2+3\delta+o(1)}~.
$$
Therefore,
$$ \sum_{u \in \mathcal{C}_1^0}\sum_{v \in \mathcal{C}_1^0} \mathrm{Cov}(\mathbf{1}_{B_u^T},\mathbf{1}_{B_v^T})
\leq n^{2(1-\delta)} n^{-2+3\delta+o(1)} = n^{\delta+o(1)} =
o\left((\mathbb{E}Y)^2\right)~.
$$
As $\mathrm{E}Y = \omega(1)$, applying Chebyshev's inequality gives:
$$\Pr[Y=0] \leq \frac{\mathrm{Var}Y}{(\mathbb{E}Y)^2} \leq \frac{\mathbb{E}Y + \sum_u\sum_v
\mathrm{Cov}(\mathbf{1}_{B_u^T},\mathbf{1}_{B_v^T})}{(\mathbb{E}Y)^2}
= \frac{1}{\mathbb{E}Y} + o(1) = o(1)~.$$ This completes the proof
of the claim and of Lemma \ref{lem-conn-or-model}.
\end{proof}

\subsection{Proof of Lemma \ref{lem-conn-and-model}}
Let $0 < \epsilon < 1$. The bounds $\tcand \leq (1+\epsilon)
\max\{\frac{1}{2},K\}$ and $\tcand \geq K - \epsilon$ will follow
from arguments similar to the ones in the proof of Lemma
\ref{lem-conn-or-model}, whereas the bound $\tcand \geq
\frac{1}{2}-\epsilon$ requires more work.

To prove that $\tcand \leq (1+\epsilon)\max\{\frac{1}{2},K\}$, take
$\epsilon'
> 0$ which satisfies:
\begin{equation}\label{eq-and-epsilon'}
\left\{\begin{array}{ll}
\ds{(1-\epsilon') \geq \frac{\left(1+\epsilon/4\right)^2}{1+\epsilon}}  \\
 \ds{1 - (1-\epsilon')^2 + (1-\epsilon')^2 K < \left(1 +
\frac{\epsilon}{4} \right) K}
\end{array}\right.~.\end{equation}
For instance, the following choice is suitable: $$ \epsilon' <
\min\left\{\frac{\epsilon K}{4}, 1-\sqrt{1-\frac{\epsilon
}{2}\min\{K,1\}}\right\} ~.$$ As before, take $c=c(\epsilon')$ such
that, with high probability, $\Gorig(c)$ has a giant component of
size at least $(1-\epsilon')n$. Defining $H^i$ to be the largest
component of $\Gand^i$, Theorem \ref{thm-M-bounded} implies that
$|H^{c'n/2}|\geq(1-\epsilon')n$ with high probability, where
$c'=\lceil\max\{K,\frac{1}{K}\}\rceil c$, and therefore we may
assume that this indeed holds. Let $A_u^i$ denote the event that a
vertex $u \in V \setminus H^i$, which belongs to some connected
component $C$, joins the giant component at time $i+1$. The
following holds:
\begin{equation}\label{eq-and-a-C} \Pr[A_u^i]=\frac{|H^i|w_C
}{\binom{n}{2} + (K-1)\binom{n-|\mathcal{C}_1^i|}{2}-K i} \mbox{ ,
where }w_C=\left\{\begin{array}{cl}
             1 & \mbox{if }|C|=1 \\
             K |C| & \mbox{otherwise}
           \end{array}\right.
~.\end{equation} By the choice of $c'$, for all $i\geq c'n/2$ we
get:
$$\Pr[A_u^i] \geq \frac{w_C (1-\epsilon')n}{\binom{n}{2}+(K-1)\binom{n-|\mathcal{C}_1^i|}{2}}
\geq \frac{2 w_C (1-\epsilon')}{n
\left(1-(1-\frac{|\mathcal{C}_1^i|}{n})^2+K(1-\frac{|\mathcal{C}_1^i|}{n})^2\right)}~.$$
For $K \geq 1$, the denominator in the expression above is clearly
bounded from above by $n K$, and for $0 < K < 1$ it is bounded from
above by $n\left(1-(1-\epsilon')^2+(1-\epsilon')^2 K\right)$, hence
\eqref{eq-and-epsilon'} implies:
\begin{equation}\label{eq-and-a-u}
\Pr[A_u^i] \geq \frac{2w_C(1-\epsilon')}{n K
(1+\epsilon/4)}\geq\frac{2w_C(1+\epsilon/4)}{n K (1+\epsilon)}~.
\end{equation}
Take $T_0=c'n/2$ and
$T=\max\{\frac{1}{2},K\}(1+\epsilon)\frac{n}{2}\log n$, and define
the following event for every $u \notin H^{T_0}$:
$$B_u^T = (u \notin H^T)~.$$ The
definition of $w_C$ in \eqref{eq-and-a-C} implies that $w_C \geq
\min\{1, 2K\}$ and when combined with \eqref{eq-and-a-u} this gives:
$$ \Pr[B_u^T] \leq \prod_{i=T_0}^T
\left(1 - \frac{2\min\{1,2K\}(1+\epsilon/4)}{n K(1+\epsilon)}
\right)\leq
$$
$$ \leq
\exp\left(-\frac{1+\epsilon/4}{1+\epsilon}\cdot\frac{\min\{\frac{1}{K},2\}}{n/2}
\left(\max\{\frac{1}{2},K\}(1+\epsilon)\frac{n}{2}\log n -
c'\frac{n}{2}\right)\right)=n^{-(1+\epsilon/4)+o(1)}~.$$ Therefore,
the expected number of vertices of $V \setminus H^{T_0}$ which do
not join $H^T$ is clearly $o(1)$, implying that
$\Gand[K][(1+\epsilon)\max\{\frac{1}{2},K\}]$ is almost surely
connected.

The lower bound $\tcand \geq K - \epsilon$ follows from the
following claim, analogous to Claim \ref{clm-iso-exponent-or}:
\begin{claim}\label{clm-iso-exponent-and} Let $0 < \delta < \frac{1}{2}$, and let
$\{H_0\}_n$ denote a family of arbitrary graphs on $n$ vertices with
$(1+o(1))n^{1-\delta}$ isolated vertices. If $\mathcal{H} \sim
\Gand|_{H_0}$, then $\mathcal{H}[K-2\delta]$ almost surely contains
isolated vertices.
\end{claim}
\begin{proof} As in the proof of Claim
\ref{clm-iso-exponent-or}, let $B_u^i$ denote the event ($u$ is
isolated at time $i$) and let $B_{u,v}^i=B_u^i \wedge B_v^i$, where
$u,v\in V(H_0)$. Clearly:
$$\Pr[\neg B_u^{i+1} | B_u^i]
= \frac{ n-1 }{\binom{n}{2} + (K-1)\binom{n-|\mathcal{C}_1^i|}{2}-K
i} ~.$$ The assumption on $H_0$ gives $|\mathcal{C}_1^i| \leq
|\mathcal{C}_1^0| = (1+o(1))n^{1-\delta}$, and thus for all $i =
O(n\log n)$:
\begin{equation}\label{eq-and-pr-b-u}\Pr[\neg B_u^{i+1} | B_u^i] =
\frac{2}{K n\left(1+O(\frac{\log n}{n}) + O(n^{-\delta})\right)}=
\frac{2+O(n^{-\delta})}{K n}~,\end{equation} and:
\begin{equation}\label{eq-and-pr-b-u-v}
\Pr[\neg B_{u,v}^{i+1} | B_{u,v}^i] = \frac{ 2n-3}{\binom{n}{2} +
(K-1)\binom{n-|\mathcal{C}_1^i|}{2}-K i} = \frac{4+O(n^{-\delta})}{K
n}~.\end{equation} Henceforth, a similar calculation to the one in
the proof of Claim \ref{clm-iso-exponent-or} gives the required
result.
\end{proof}
Choosing $\delta=\epsilon/2$ and applying the last claim implies
that $\tcand \geq K - \epsilon$, and it remains to show that $\tcand
\geq \frac{1}{2}-\epsilon$. This follows from the fact that
$\Gand[K][\frac{1}{2}-\epsilon]$ contains components of size $2$
almost surely. Unfortunately, we cannot repeat the last argument in
order to show that indeed this is the case, since we have no
guarantee that at time $\tau$, the hitting time for the property
$|\mathcal{C}_1|\leq n^{1-\delta}$, $\Gand$ still satisfies
$|\mathcal{C}_2|=n^{1-\alpha}$ for some small $\alpha$. To prove the
next claim, which completes the proof of Lemma
\ref{lem-conn-or-model}, we consider a simplified process, where it
is possible to give a lower bound on $|\mathcal{C}_2|$, and show
that this process stochastically dominates $\Gand$ with regards to
$|\mathcal{C}_2|$.

\begin{claim}\label{clm-iso-exponent-and-edges} For
all $K\in(0,1]$ and $\epsilon > 0$, $\Gand[K][\frac{1}{2}-\epsilon]$
almost surely contains a connected component of size $2$.
\end{claim}
\begin{proof}[Proof of claim]
Define the process $\Gandap$, which is an {\em approximated} version
of $\Gand$, as follows: at each step, assign a weight $K$ to {\em
ordered pairs} of the form $\{(u,v):u,v \notin \mathcal{C}_1\}$, and
a weight $1$ otherwise. If the ordered pair $(u,v)$ chosen at some
step $i$ is a (self) loop or corresponds to an edge which already
exists in $\Gandap$, this step is omitted. In other words, $\Gandap$
assigns weights to all ordered pairs, and disregards selections of
multiple edges or loops. Clearly, if $\Gandap^t$ contains $m \leq t$
edges, then $\Gandap^t \sim \Gand^m$. Furthermore, for any
$t=o(n^2)$ we have:
$$ \mathbb{E}\left(t-|E(\Gandap^t)|\right) \leq \sum_{i=1}^t
\max\{\frac{1}{K},K\} \frac{t + n}{n^2} =
\left(O\left(\frac{t}{n^2}\right)+O\left(\frac{1}{n}\right)\right)t
= o(t)~,$$ where $E(H)$ denotes the set of edges of the graph $H$.
In particular, for any fixed $c$, with high probability $\Gand[K](c)
\sim \Gandap[K](c+o(1))$ and $\Gand[K][c] \sim \Gandap[K][c+o(1)]$,
and it is sufficient to prove the claim for $\Gandap$.

Take $0 < \epsilon < \frac{1}{2}$. By well known results on the
original Erd\H{o}s-R\'enyi graph process, for any fixed $t$, with
high probability: $|\mathcal{C}_2(\Gorig(t))|=\Theta(n)$, with the
constant tending to $0$ as $t\to\infty$. For a short proof of this
fact, one can verify that the following differential equation
approximates the graph parameter $I_2$ along $\Gorig$ up to an
$o(1)$ error term:
$$ w'(t) = -2w(t) + y^2(t)~,~w(0)=0~,$$ where $y(t)$ is the approximating
function of the fraction of isolated vertices, $I(t)$, and the
time-scaling is of $n/2$ edges at each step. Substituting the well
known fact (which also follows from Theorem \ref{thm-tgand} when
substituting $K=1$) that $y(t)=\exp(-t)$, it follows that
$w(t)=t\exp(-2t)$.

By Theorem \ref{thm-M-bounded} and the above fact, we deduce that
for every sufficiently large $c$ there exist $0 < \alpha_1 <
\alpha_2 < \frac{1}{3}$ such that $\Gandap(c)$ almost surely
satisfies $\alpha_1 n \leq |\mathcal{C}_2| \leq \alpha_2 n$. Let:
$$c' = \inf\{t > c : y(t) \leq \frac{\epsilon K}{2}\}~, $$ where $y$
is the solution to the ODE \eqref{y-eq}), and set $\epsilon' =
y(c')$. By Theorem \ref{thm-tgand}, with high probability
$I\left(\Gandap(c')\right) = y(c') + o(1)$, and in particular,
$\Gandap(c')$ has $(1+o(1))\epsilon' n$ isolated vertices almost
surely. Let $H_0$ denote $\Gandap(c')$, and take
$\alpha_1=\alpha_1(c')$, $\alpha_2=\alpha_2(c')$ as above. According
to these definitions, $S=\mathcal{C}_1(H_0)$ satisfies $s=|S| =
(1+o(1))\epsilon' n$ almost surely, and $W = \mathcal{C}_2(H_0)$
satisfies $\alpha_1 n \leq |W| \leq \alpha_2 n $ almost surely.
Assume that indeed this holds.

We consider the graph process $\mathcal{H}$, which begins with
$H_0$, and at each step selects an ordered pair $(u,v) \in V(H_0)^2$
according to the following probabilities:
$$ \left\{\begin{array}{ll}
 \frac{1}{K n^2} & \mbox{if $(u,v)$ is incident to $W$ and $S$}
 \\
 \frac{1}{n^2} & \mbox{if $(u,v)$ is incident to $W$ and $V \setminus
 S$} \\
 \lambda & \mbox{otherwise}
\end{array}\right.~,$$
where the value of $\lambda>0$ is chosen such that the probabilities
sum up to $1$. This is possible, since the probabilities for the
first two types of pairs sum up to at most:
 $$2\alpha_2 \left(\epsilon'/K + (1-\epsilon')\right) \leq 3 \alpha_2 < 1~, $$
as $\epsilon'\leq\frac{K}{2}$ and $\alpha_2<\frac{1}{3}$.
 We claim that the following two events
occur almost surely, and complete the proof of the claim:
\begin{eqnarray}
&\mathcal{C}_2(\mathcal{H}[\frac{1}{2}-\epsilon]) \cap W \neq \emptyset~, \label{eq-h-has-iso-edges}\\
&\left|\mathcal{C}_2(\mathcal{H}[\frac{1}{2}-\epsilon]) \cap
W\right| \leq
\left|\mathcal{C}_2(\Gandap|_{H_0}[\frac{1}{2}-\epsilon]) \cap
W\right|~. \label{eq-h-contained-g}
\end{eqnarray}

A standard second moment consideration proves that
\eqref{eq-h-has-iso-edges} occurs with high probability. To see
this, set: $T = \left(\frac{1}{2}-\epsilon\right)\frac{n}{2}\log n$,
and notice that:
\begin{equation}\frac{2T}{Kn^2}\left(K(n-s-2)+s\right)=
\frac{1}{2}(1-2\epsilon)\left((1-\epsilon')+\frac{\epsilon'}{K}+o(1)\right)
\log n =\frac{1}{2}(1-\epsilon''+o(1))\log n
\end{equation} for some $\epsilon \leq \epsilon'' \leq 2\epsilon$ (by our
choice of $\epsilon'$). Next, let the random variable $X_e$ ($e\in
W$) be the indicator of the event: $\left(e\in
\mathcal{C}_2(\mathcal{H}[\frac{1}{2}-\epsilon])\right)$, and let
$X=\sum_{e\in W}X_e$. The following holds:
$$ \Pr[X_e=1] = \left(1- \frac{2K(n-s-2)+2s}{Kn^2}\right)^T ~,$$
$$ \Pr[X_e=1 \wedge X_{e'}=1] = \left(1- \frac{4K(n-s-2)-4K+4s}{Kn^2} \right)^T ~,$$
for every $e,e' \in W$. Thus, a calculation similar to the one in
the proof of Claim \ref{clm-iso-exponent-or} gives:
$$ \mathbb{E}X \geq |W| \exp\left(-(2-o(1))T\frac{K(n-s-2)+s}{Kn^2}\right) \geq
|W| n^{(-1+\epsilon'')/2+o(1)} = \omega(\sqrt{n})~,$$ and:
$$ \mathrm{Cov}(X_e, X_{e'}) \leq
\frac{4K}{Kn^2} \exp\left(-4\frac{K(n-s-4)+K+s}{Kn^2} T\right)\leq
4n^{-3+\epsilon''+o(1)}~.$$ Therefore:
$$ \sum_{e \in W} \sum_{e'\in W} \mathrm{Cov}(X_e,X_{e'}) \leq
4 |W|^2 n^{-3+\epsilon''+o(1)} = o(\mathbb{E}X)~,$$ and in
particular, $\mathrm{Var}(X)=(1+o(1))\mathbb{E}X =
o\left(\mathbb{E}X)^2\right)$ and by Chebyshev's inequality we
deduce that $X>0$ almost surely.

It remains to prove that \eqref{eq-h-contained-g} occurs almost
surely. This is achieved by a coupling argument: we claim that there
exists a coupling of the processes $\Gandap|_{H_0}$ and
$\mathcal{H}$ whose support consists of pairs $(G_t,H_t)$ such that:
$G_t \sim \Gandap|_{H_0}^t$, $H_t \sim \mathcal{H}^t$, and
$\left(\mathcal{C}_2(H_t) \cap W\right) \subset
\left(\mathcal{C}_2(G_t) \cap W\right)$. This clearly holds for
$t=0$, and by induction, it remains to extend the coupling from
$(G_t,H_t)$ to $(G_{t+1},H_{t+1})$. For this purpose, we apply the
following lemma of \cite{AGLS} (Lemma 2.2), which was first proved
by Strassen \cite{Strassen} in a slightly different setting:
\begin{lemma}[\cite{AGLS}]\label{lemma-coupling}
Let $U,V$ be two finite sets, and let $R\subset U \times V$ denote a
relation on $U,V$. Let $\mu$ and $\nu$ denote probability measures
on $U$ and $V$ respectively, such that the following inequality
holds for every $A \subset U$:
\begin{equation}\label{eq-hall-condition}
\mu(A) \leq \nu(\{y\in V : x R y \mbox{ for some }x \in A\})~.
\end{equation}
Then there exists a coupling $\varphi$ of $\mu,\nu$ whose support is
contained in $R$.
\end{lemma}
Let $U,V$ be the set of all $n^2$ ordered pairs selected,
representing the next pair selected by $G_t$ and $H_t$ respectively.
Let $\mu$ denote the probability measure of each selection in $U$ by
$G_t$, and let $\nu$ denote the probability of each selection in $V$
by $H_t$. Define: $$ \left\{\begin{array}{lcl}
X=\mathcal{C}_2(G_t) \cap W &,& X_{u,v}=\mathcal{C}_2\left(G_t \cup (u,v)\right) \cap W\\
Y=\mathcal{C}_2(H_t) \cap W &,& Y_{x,y}=\mathcal{C}_2\left(H_t \cup
(x,y)\right) \cap W \end{array}\right.~.$$ By the induction
hypothesis, $Y \subseteq X$, and we define $R$ to be $\{
\left((u,v),(x,y)\right) : Y_{x,y} \subseteq X_{u,v} \}$.

Clearly, if $(u,v)$ is a loop or an edge which already belongs to
$G_t$, it has no effect on $X$, and as $Y_{x,y} \subseteq Y
\subseteq X$ we get $(u,v)R(x,y)$ for all $x,y$. Furthermore, every
$(u,v)$ which is not incident to any $e \in Y$ also satisfies
$(u,v)R(x,y)$ for all $x,y$, as the components $(u,v)$ may remove
from $X$ already do not belong to $Y$. Therefore,
\eqref{eq-hall-condition} is satisfied for every $A \subset U$ which
contains such pairs $(u,v)$.

It remains to prove that \eqref{eq-hall-condition} holds for sets $A
\subset U$ such that $A \cap E(G_t) = \emptyset$, and $A$ consists
entirely of edges incident to edges of $Y$. Let
$A=\{e_1,\ldots,e_m\}$, and notice that $e_i \notin E(H_t)$ for all
$i$, as $e_i$ is incident to some component $e \in Y$ of size $2$ in
$H_t$, satisfying $e\neq e_i$ as $Y \subset E(G_t)$. Furthermore,
for all $i$ we have $e_i R e_i$, since all components that $e_i$
removes from $X$ in $X_{e_i}$ are also removed from $Y$ in
$Y_{e_i}$. Thus, showing that $\mu(e_i) \leq \nu(e_i)$ for all $i$
will imply that $A$ satisfies the condition of
\eqref{eq-hall-condition}. Indeed, if $e_i=(u,v)$ is such that $u,v
\notin S$ ($u,v$ are both non-isolated in $H_0$), then:
$$ \mu(e_i) = \frac{K}{n^2 +
(K-1)(n-\mathcal{C}_1(G^t))^ 2} \leq \frac{K}{K n^2} = \nu(e_i)~,$$
and otherwise:
$$ \mu(e_i) \leq \frac{1}{n^2 +
(K-1)(n-\mathcal{C}_1(G^t))^ 2} \leq \frac{1}{K n^2} = \nu(e_i)~,$$
by definition of the process $\mathcal{H}$, completing the proof of
the claim and the proof of Lemma \ref{lem-conn-or-model}.
\end{proof}

\subsection{Proof of Lemma \ref{lem-conn-K-0}}
To prove that $\tcor(0)=\frac{1}{2}$, we recall the following easy
facts stated in \cite{AGLS}: by definition, $\Gor[0]$ adds edges
between pairs of isolated vertices until no such pair is left.
Hence, after adding $\lfloor n/2\rfloor$ edges there is at most $1$
isolated vertex in $\Gor$, and $\Gor$ behaves as $\Gorig$ on
$\lfloor n/2\rfloor$ components of size 2 (and possibly $1$
additional isolated vertex). Thus,
$\tgor(0)=1+\frac{1}{2}=\frac{3}{2}$.

Proving a connectivity threshold of
$(\frac{1}{2}+o(1))\frac{n}{2}\log n$, we may assume $n$ is even:
otherwise, the single isolated vertex from time $\lfloor n/2
\rfloor$ becomes connected almost surely at time $\omega(n)$, and
its edge set accounts to at most $n-1 = o(n \log n)$ edges, not
affecting the threshold for connectivity. By the discussion above,
for even values of $n$, $\Gor[0]$ has a linear number of components
of size $2$ and no isolated vertices at the time of appearance of
the giant component. Thus, we deduce from the arguments used for the
proofs of Lemmas \ref{lem-conn-or-model} and
\ref{lem-conn-and-model} that, with high probability,
$\Gor[0][\frac{1}{2}-\epsilon]$ still contains components of size
$2$, whereas $\Gor[0][\frac{1}{2}+\epsilon]$ is connected.
Altogether, $\tcor(0)=\frac{1}{2}$.

It remains to show that $\tcand(0)=\frac{1}{2}$. Substituting $K=0$
in equation \eqref{y-eq}, it reduces to the form:
$$y'=\frac{1}{y-2}~,~y(0)=1~,$$ provided that $y \neq 0$. This yields the
solution:
\begin{equation}\label{eq-y-sol-K=0}y(t)=2-\sqrt{1+2t}~.\end{equation}
Notice that $y(t)$ is strictly monotone decreasing, and reaches $0$
at $t=\frac{3}{2}$. Therefore, if we denote by $\tau_0(\Gand[0])$
the minimal time $t$ at which $I\left(\Gand[0](t)\right)=0$, Theorem
\ref{thm-tgand} implies that for every $0 < \delta <  1$, $|\tau -
\frac{3}{2}| < \delta$ almost surely. To see this, let $x_0$ be such
that $y(x_0)=\delta$, and let $x_1 = \max\{ x_0,
\frac{3}{2}-\delta\}$. By Theorem \ref{thm-tgand}, with high
probability:
$$I(\Gand[0](\frac{3}{2}-\delta)) \geq I\left(\Gand[0](x_1)\right) = y(x_1)+o(1) \geq
\frac{1}{2} y(x_1) > 0~,$$ where the second from last inequality
holds for sufficiently large values of $n$. On the other hand,
$y(x_1) < \delta$, and hence
$I\left(\Gand[0](\frac{3}{2}+\delta)\right) \leq
I\left(\Gand[0](x_1+\delta)\right) = 0$ almost surely, since each
edge eliminates at least one isolated vertex in $\Gand[0]$.

Equation \eqref{w-eq} takes the following form after substituting
$K=0$ and the value of $y(t)$ as it is given in
\eqref{eq-y-sol-K=0}, provided that $y \neq 0$:
$$ w'(t) = \frac{y-2w}{2-y}=2\frac{1-w}{\sqrt{1+2t}}-1~,~w(0)=0~,$$
which has the solution: \begin{equation}\label{eq-w-sol-K=0}
w(t)=\frac{5}{4}-\frac{3}{4}\mathrm{e}^{2(1-\sqrt{1+2t})}
-\frac{1}{2}\sqrt{1+2t}~.
\end{equation}
For $t=\frac{3}{2}$ we have
$w(t)=\frac{1}{4}-\frac{3}{4\mathrm{e}^2}\approx 0.1485$, hence from
the discussion above we obtain that $H=\Gand(\tau)$ almost surely
satisfies $I_2(H)=\Theta(n)$ (and $I(H)=0$). From that point on, the
process $\Gand$ is equivalent to $\Gorig|_H$, and from the similar
arguments to those used in the proofs of Lemmas
\ref{lem-conn-or-model} and \ref{lem-conn-and-model} we obtain that
$\tcor(0)=\frac{1}{2}$. \qed

\subsection{Computer experiments of $\mathbf{\tcand}$ and $\mathbf{\tcor}$}\label{sec::sim1}
Maintaining the set of isolated vertices and the edges already added
to the process provides all the information needed to add the next
edge to $\Gand$ and $\Gor$ at an $O(1)$ cost. In order to recognize
the threshold for connectivity, the set of connected components must
be efficiently maintained. Our implementation holds the components
in linked-lists, according to the Weighted-Union Heuristic (see,
e.g., \cite{CLR} p. 445) which guarantees an average cost of $O(\log
n)$ for uniting components.

\begin{figure}
\centering
\includegraphics{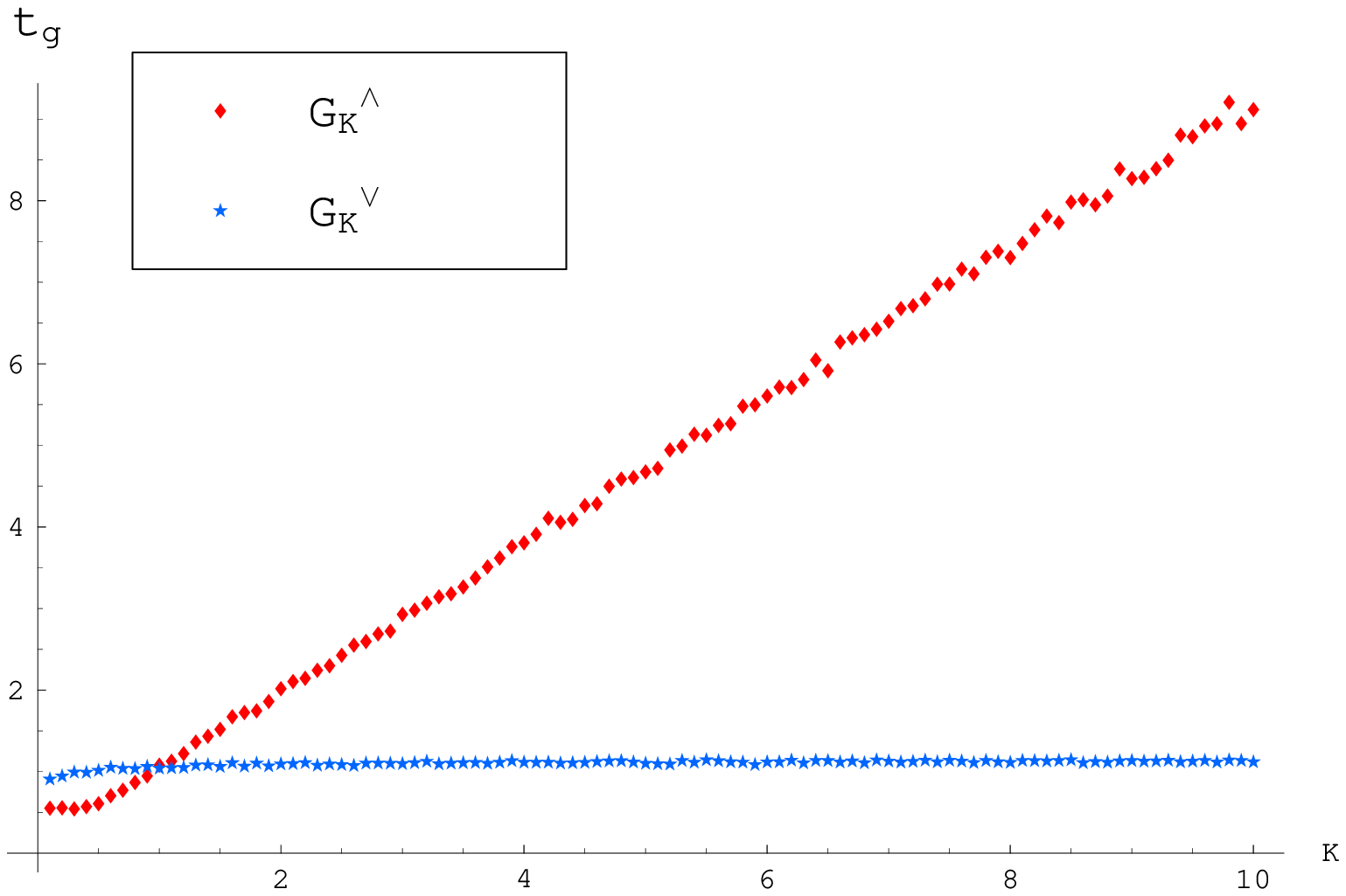}
\caption{Comparison of the numerical results for $\tcand$ and
$\tcor$, and estimations of $\tcand$ and $\tcor$ according to
computer simulations of the model.} \label{fig::sim-conn}
\end{figure}

Figure \ref{fig::sim-conn} shows the results of $\tcand$ and $\tcor$
according to simulations of both models on $n=10^4$ vertices. The
values of $\tcand$ and $\tcor$ were averaged over 100 tests per
value of $K$.

\section{The appearance of a giant component in $\mathbf{\Gand}$ }\label{sec::giant}
\subsection{Proof of Theorem \ref{thm-tgand}}
We begin with a short summary of the methods used in \cite{AGLS} to
analyze $\tgor(K)$ for $K>0$. First, the authors prove Theorem
\ref{thm-M-bounded} and deduce that the process $\Gor$ is
stochastically dominated by $\Gorig$, up to a timescale factor of
$\lceil \max\{\frac{1}{K},K\} \rceil$. By applying the differential
equation method of Wormald \cite{WormaldDiffEq} to the approximated
process $\Gorap$ (which selects an ordered pair at each step,
similar to the process $\Gandap$ introduced in the proof of Claim
\ref{clm-iso-exponent-and-edges}), the parameters $I(G)$ and $S(G)$
are approximated by $y(t)$ and $z(t)$, solutions to a system of
coupled ODEs. Henceforth, a repeated use of Theorem 3.1 of
\cite{SpencerWormald}, which relates the susceptibility and the
appearance of a giant component, implies that the singularity point
$z(t)$ is equal to $\tcor$.

We note that the methods of \cite{AGLS} can be applied to any family
of $M$-bounded processes $\mathcal{P}_M$ (as referred to in Theorem
\ref{thm-M-bounded}) provided that the following conditions hold:
\begin{enumerate}
  \item Attempting to approximate a fixed number of bounded graph parameters
  $Y_1,\ldots,Y_d$ (such as $I(G)$ and $I_2(G)$) by functions
  $y_1,\ldots,y_d$ for $1 \leq t \leq T$, we require that
  \begin{equation}\label{eq-diff-eq-requirement}
  \mathbb{E}\frac{Y_i(\mathcal{P}_M^{t+1})-Y_i(\mathcal{P}_M^t)}{2/n} =
  \frac{d y_i}{d t}\left(Y_1(\mathcal{P}_M^t),\ldots,Y_d(\mathcal{P}_M^t),t\right) + \mathrm{err_{y_i}(t)}~,
  \end{equation} for all $i$ and $t$, where $\max_{i,t} \mathrm{err}_{y_i}(t)=o(1)$.

  \item In addition to the above approximations of $Y_1,\ldots,Y_d$ by $y_1,\ldots,y_d$,
  when attempting to approximate $S(G)$ by the function $z(t)$ for $1 \leq t \leq T$, we require that
  \begin{equation}\label{eq-diff-eq-requirement-z}
  \mathbb{E}\frac{S(\mathcal{P}_M^{t+1})-S(\mathcal{P}_M^t)}{2/n} =
  \frac{d z}{d t}\left(Y_1(\mathcal{P}_M^t),\ldots,Y_d(\mathcal{P}_M^t),S(\mathcal{P}_M^t),t\right) + \mathrm{err_{z}(t)}~,
  \end{equation} where $\max_{t} \mathrm{err}_{z}(t)=o(1)$.

  To show the above, we may use the fact that with high probability, the largest component is of size $O(\log n)$ as long as $t <
x_c -\epsilon$, where $x_c$ is a possible singularity point of
$z(t)$. This follows from the proof of Theorem 1.3 of \cite{AGLS},
which applies to this generalized setting as well.

  \item Finally, if the above function $z(t)$ has a singularity point, we require that it is uniformly bounded (regardless of $M$).
\end{enumerate}
If the above 3 conditions hold, we obtain that $Y_1,\ldots,Y_d,S$
are within $o(1)$ distance from $y_1,\ldots,y_d,z$ along the process
$\mathcal{P}_M$ for $1\leq t \leq T$. Furthermore, if $z$ has a
singularity point at $x_c$ and the above conditions hold with
$T=x_c-\epsilon$ for any $\epsilon>0$, it follows that the
appearance of a giant component in $\mathcal{P}_M$ is at $t=x_c$.

We first show that the above 3 conditions hold for $K>0$ under the
assumptions of Theorem \ref{thm-tgand}. By the well known fact that
for any constant $t$, $I(\Gorig(t))=\Theta(n)$ almost surely,
Theorem \ref{thm-M-bounded} implies that $\tau_\delta$ can be taken
to be arbitrarily large. In particular, we can take $\tau_\delta >
x_c$. In order to verify that \eqref{eq-diff-eq-requirement} holds
for $I(G)$ and $I_2(G)$, let $G$ and $G'$ denote $\Gand^t$ and
$\Gand^{t+1}$ respectively for some fixed $t$, and let $I=I(G)$ and
$I_2=I_2(G)$. We have:
\begin{eqnarray}\mathbb{E}\frac{I(G')-I(G)}{2/n} &=&
(-1)\frac{I(I-\frac{1}{n})}{1+(K-1)(1-I)^2}+
(-1)\frac{2I(1-I)}{1+(K-1)(1-I)^2} = \nonumber \\
&=& \frac{-I}{1+(K-1)(1-I)^2} + O(I/n) \label{eq-I-diff-eq-delta}~,
\end{eqnarray}
and:
\begin{eqnarray}\mathbb{E}\frac{I_2(G')-I_2(G)}{2/n} &=& \frac
{I(I-\frac{1}{n}) - 2I_2 I - 2K I_2(1-I-I_2) - 2 K
I_2(I_2-\frac{2}{n})}{1+(K-1)(1-I)^2}
= \nonumber \\
&=& \frac{I^2 -2I_2 I -2K I_2(1-I)}{1+(K-1)(1-I)^2} +
O(I/n)+O(I_2/n)) \label{eq-I2-diff-eq-delta}~,
\end{eqnarray}
where in both cases we used the fact that $K>0$ to obtain an upper
bound of $o(1)$ on the error term. To prove
\eqref{eq-diff-eq-requirement-z}, set $S=S(G)$ and observe that:
$$ \mathbb{E}\frac{S(G')-S(G)}{2/n} = $$
\begin{eqnarray}
&=&\frac{n/2}{1+(K-1)(1-I)^2} \bigg(I(I-\frac{1}{n})\frac{2}{n} +
 2I \mathop{\sum_{C\in\mathcal{C}}}_{|C|>1}
\frac{|C|}{n}\frac{2|C|}{n} +
\mathop{\sum_{C_1\in\mathcal{C}}}_{|C_1|>1}
\mathop{\sum_{C_2\in\mathcal{C}\setminus\{C_1\}}}_{|C_2|>1} \frac{K
|C_1||C_2|}{n^2}\frac{2|C_1||C_2|}{n}\bigg)= \nonumber \\
&=&\frac{2I(S-I)+I^2+K(S-I)^2}{1+(K-1)(1-I)^2}+O(I/n)+O\bigg(\mathop{\sum_{C\in\mathcal{C}}}_{|C|>1}\frac{|C|^4}{n^2}\bigg)=\nonumber\\
&=&\frac{S^2+(K-1)(S-I)^2}{1+(K-1)(1-I)^2}+O(I/n)+O\bigg(\mathop{\sum_{C\in\mathcal{C}}}_{|C|>1}\frac{|C|^4}{n^2}\bigg)\label{eq-S-diff-eq-delta}~,
\end{eqnarray}
and the assumption that $|C|=O(\log n)$ gives a bound of $o(1)$ on
$\mathrm{err}_z$. The next claim therefore completes the statements
of Theorem \ref{thm-tgand} for the case $K>0$:

\begin{claim}For every $K > 0$, $\tgand(K) < 5$.
\end{claim}\begin{proof}
Consider the case $K \geq 1$, and let $y(t)$ and $z(t)$ denote the
solutions to the ODEs \eqref{y-eq} and \eqref{z-eq} respectively.
Recalling that $y(t) \leq 1$ for every $t\geq 0$, \eqref{z-eq}
yields that $z'(t) \geq 0$ provided that $z(t)\geq 1$. Thus, the
initial condition $z(0)=1$ implies that $z(t)$ is monotone
increasing in $t$, and in particular:
\begin{equation}\label{eq-z-K-geq-1}z' =
\frac{z^2+(K-1)(z-y)^2}{1+(K-1)(1-y)^2} \geq
\frac{z^2+(K-1)z^2(1-y)^2}{1+(K-1)(1-y)^2} = z^2~,\end{equation}
where inequality is by the fact that $y \geq 0$ and $z \geq 1$. By
standard considerations in differential analysis ,
\eqref{eq-z-K-geq-1} and the initial condition $z(0)=1$ imply that
$z(t) \geq \frac{1}{1-t}$ for every $t\geq 0$ (as
$\tilde{z}(t)=\frac{1}{1-t}$ satisfies $\tilde{z}'=\tilde{z}^2$ and
$\tilde{z}(0)=1$), and in particular, $\tgand \leq 1$ for any $K
\geq 1$.

We are left with the case $0 < K < 1$. Clearly, \eqref{y-eq} implies
that $y' \leq 0$ for every $t \geq 0$, and furthermore, $y' < 0$ as
long as $y>0$, hence $y$ is strictly monotone decreasing from $1$ to
$0$. Let $t^*$ be such that
\begin{equation}\label{eq-y-t*}
y(t^*)=1-\sqrt{1-K}~.\end{equation} Notice that the solutions to the
equation $-x^2+2x-\frac{K}{1-K} = 0$ are: $x_{1,2}=1\pm
\sqrt{1-\frac{K}{1-K}}$ if $0 < K \leq \frac{1}{2}$, and no solution
exists if $\frac{1}{2} < K < 1$. In both cases,
$-x^2+2x-\frac{K}{1-K} < 0$ for every $x \leq 1-\sqrt{1-K}$.
Therefore, the definition of $t^*$ and the fact that $y$ is monotone
decreasing give:
$$ -y(t)^2 +2y(t) - \frac{K}{1-K} \leq 0 ~\mbox{ for every }t \geq t^*,$$ or equivalently:
\begin{equation}\label{eq-y-bound-t*}(1-K)y(t)(2-y(t)) \leq K ~\mbox{ for every }t\geq t^*~.\end{equation}
Rewriting equation \eqref{z-eq} as:
\begin{equation}\label{eq-z-K-leq-1}
z' = \frac{z^2 - (1-K)(z-y)^2}{1-(1-K)(1-y)^2} = \frac{K z^2 +
(1-K)y(2z-y)}{K + (1-K)y(2-y)}~,
\end{equation}
it follows that for every $t \geq t^*$, $z' \geq \frac{1}{2}z^2$.
Furthermore, for every $t \geq 0$, $z' > 0$, and hence $z(t^*)
> z(0) = 1$. We obtain that the function $w(t) = z(t-t^*)$ satisfies
$w' \geq \frac{1}{2}w^2$ for every $t \geq 0$ and $w(0)\geq 1$, and
by the same consideration as above, $w(t) \geq \frac{2}{2-t}$ for
every $t \geq 0$. Altogether, we deduce that $\tgand \leq t^* + 2$,
and it remains to provide an upper bound on $t^*$.

For this purpose, define $u(t)=1-y(t)$, and consider \eqref{y-eq}
for $0 \leq t \leq t^*$:
$$ u' = -y' = \frac{1-u}{1-(1-K)u^2} \geq \frac{1-u}{1-u^4} =
\frac{1}{1+u+u^2+u^3}~,$$ where the inequality is by the fact that
$u(t) \leq \sqrt{1-K}$ for $0\leq t \leq t^*$. Define $w(t)$ to be
the solution to the differential equation:
\begin{equation}\label{eq-w-cubic}w' = \frac{1}{1+w+w^2+w^3}~,~w(0)=0
~,\end{equation} it follows from the above mentioned argument that
$u(t) \geq w(t)$ for $0\leq t \leq t^*$. The solution to
\eqref{eq-w-cubic} satisfies: $t = \sum_{j=1}^4 \frac{w^j}{j}$,
hence $w(t_0) = \sqrt{1-K}$ for $t_0 = \sum_{j=1}^4
\frac{(1-K)^{j/2}}{j} \leq \frac{25}{12}$. As $u(t_0) \geq w(t_0)$,
it follows that $t^* \leq t_0 \leq \frac{25}{12}$, completing the
proof.
\end{proof}

In the special case $K=0$, $\Gand[0]$ is no longer an $M$-bounded
process, however the assertions of statements $1,2$ of Theorem
\ref{thm-tgand} remain valid and follow from
\eqref{eq-diff-eq-requirement} and \eqref{eq-diff-eq-requirement-z},
by applying Wormald's differential equation method directly. To see
that \eqref{eq-diff-eq-requirement} holds, notice that as long as
$I(\Gand(t)) \geq \delta$ for some fixed $\delta > 0$, the
denominators in \eqref{eq-I-diff-eq-delta} and
\eqref{eq-I2-diff-eq-delta} are $\Theta(1)$, and the approximation
remains valid (note that $z(t)$ has no singularity point for $K=0$).
To show that \eqref{eq-diff-eq-requirement-z} holds, set $S=S(G)$,
and note that:

\begin{eqnarray}\mathbb{E}\frac{S(G')-S(G)}{2/n} &=&\frac{n/2}{1-(1-I)^2} \bigg(I(I-\frac{1}{n})\frac{2}{n} +
 2I \mathop{\sum_{C\in\mathcal{C}}}_{|C|>1}
\frac{|C|}{n}\frac{2|C|}{n} \bigg)= \nonumber \\
&=&\frac{2I(S-I)+I^2}{2I-I^2}+O(I/n)=\frac{2S-I}{2-I}+O(I/n)\label{eq-S-diff-eq-delta-K=0}~,
\end{eqnarray}
where we used the fact that $I > 0$ for $t \leq \tau_\delta$. This
implies that the error term $\mathrm{err}_z$ is $o(1)$ without
making any assumptions on the size of the largest component, and
completes the proof of the theorem. \qed

\subsection{Computer experiments of $\mathbf{\tgand}$}\label{sec::sim2}
We conducted simulations of $\tgand$ using the implementation of
$\Gand$ mentioned in \ref{sec::sim1}. In these simulations, the
number of vertices was $n=10^6$, and the threshold for the
appearance of the giant component was taken to be the minimal time
at which $\Gand$ contains a component of size $\alpha n$, where
$\alpha=0.01$. The value of $\tgand(K)$ was averaged over $10$ tests
for each value of $K$.

\begin{figure} \centering
\includegraphics{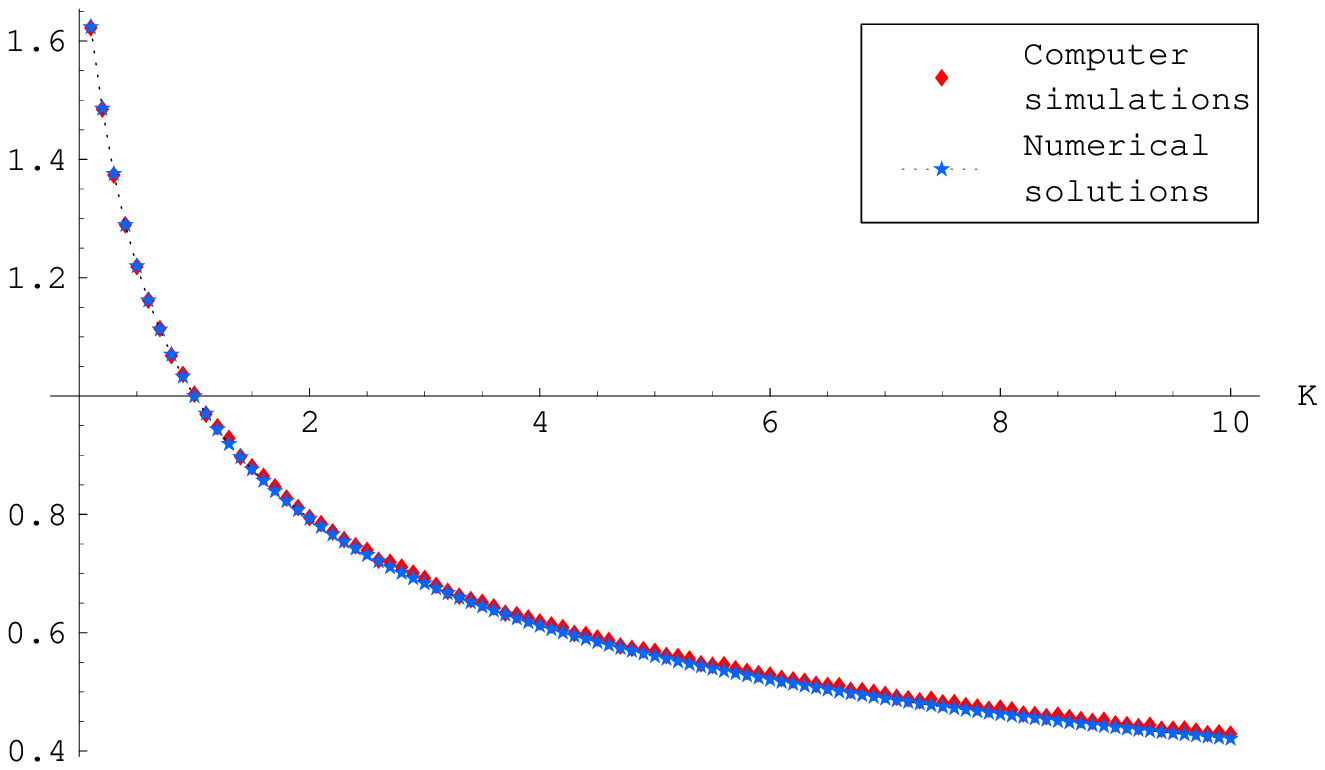}
\caption{Comparison of the numerical results for $\tgand$, and
estimations of $\tgand$ according to computer simulations of the
model.} \label{fig::sim-giant}
\end{figure}

Figure \ref{fig::sim-giant} shows the comparison between the values
of $\tgand$ according to the above computer simulations, and the
values obtained by numerically solving the ODEs \eqref{y-eq} and
\eqref{z-eq} by Mathematica.

\subsection{Proof of Theorem \ref{thm-K-0-K-gg-1}}
First, the fact that $\tcand(K)$ is continuous follows from the
general continuous dependence of ODEs on their parameters (in this
case,the single parameter $K$).

For the special case $K=0$, recall that in our treatment of
$\Gand[0]$ for the proof of Lemma \ref{lem-conn-K-0} we showed that
for every $\delta > 0$, with high probability $y(t)=2-\sqrt{1+2t}$
approximates $I(\Gand[0](t))$ for $t\leq \frac{3}{2}-\delta$, and
$I(\Gand[0](\frac{3}{2}+\delta))=0$. Take $\delta > 0$ and consider
the interval $[0,\frac{3}{2}-\delta]$; equation \eqref{z-eq} takes
the following form when substituting $K=0$ and the solution to
$y(t)$:
\begin{equation}
\label{eq-z-K=0-eq} z'
=\frac{z^2-(z-y)^2}{2y-y^2}=\frac{2z-y}{2-y}=\frac{2(z-1)}{\sqrt{1+2t}}+1~.\end{equation}
Taking $w=z-1$, we obtain the linear equation:
\begin{equation}\label{eq-linear-w}
w' - \frac{2}{\sqrt{1+2t}}w=1~. \end{equation} Multiplying
\eqref{eq-linear-w} by its integrating factor and integrating by
parts, we get:
$$ w(t) = \exp\left\{2\int (1+2t)^{-\frac{1}{2}} dt
\right\} \int \exp\left\{-2\int (1+2t)^{-\frac{1}{2}} dt
\right\}\,dx~=$$
$$ = -\frac{1}{2}(\sqrt{1+2t}+\frac{1}{2}) + C \exp(2\sqrt{1+2t})~.$$
The initial condition $w(0)=0$ gives $C=\frac{3}{4}\mathrm{e}^{-2}$,
hence:
$$z(t)=\frac{3}{4}\mathrm{e}^{2(\sqrt{1+2t}-1)}-\frac{1}{2}\sqrt{1+2t}+\frac{3}{4}
\mbox{ for all }t \in [0,\frac{3}{2}-\delta]~.$$ The above solutions
for $y(t)$ and $z(t)$ give $y(\frac{3}{2})=0$ and
\begin{equation}\label{eq-z-t0}
z(3/2)=\frac{3}{4}\mathrm{e}^2-\frac{1}{4} ~.\end{equation} Let
$\tau_0=\tau_0(\Gand[0])$ denote the first time $t$ at which
$\Gand[0]^t$ has no isolated vertices, i.e., the hitting time for
the property: $\{G: I(G)=0\}$. By the above arguments, we have:
\begin{equation}\label{eq-tau0-value}\tau_0=\left(\frac{3}{2}+o(1)\right)n/2\end{equation}
almost surely. We claim that proving that with high probability:
\begin{equation}\label{eq-susc-3/2-assumption}
S(\Gand[0]^{\tau_0})=z(3/2)+o(1) \end{equation} implies the required
result on $\tgand(0)$. Indeed, once no isolated vertices are left,
the process $\Gand[0]$ adds edges according to the uniform
distribution, and hence from that point the susceptibility follows
the equation $z'=z^2$ (e.g., see the case $K=1$ of the analysis of
$\Gand$ or $\Gor$). By \eqref{eq-tau0-value}, we obtain that for
some $x_c > \frac{3}{2}$:
$$S(\Gand[0](t)) = (1+o(1))\frac{1}{x_c-t} \mbox{ for any  }\frac{3}{2}\leq t < x_c~,$$
and the value of $x_c$ is derived from the initial condition
\eqref{eq-z-t0}: \begin{equation}\label{eq-tgand0-result}\tgand(0) =
x_c = \frac{3}{2} + \frac{4}{3\mathrm{e}^2-1}~.\end{equation} It is
left to show that \eqref{eq-susc-3/2-assumption} indeed holds. The
lower bound $S(\Gand[0]^{\tau_0})\geq z(\frac{3}{2})-o(1)$ follows
from Theorem \ref{thm-tgand}, which states that $z$ approximates $S$
until time $\frac{3}{2}-\delta$ for any $\delta>0$, and from the
continuity and monotonicity of $z$. That is, for any fixed $\xi>0$,
choosing a sufficiently small $\delta>0$ such that
$z(\frac{3}{2}-\delta)
> z(\frac{3}{2})-\xi$ gives:
$$S\left(\Gand[0](\frac{3}{2}-\delta)\right) =
z(\frac{3}{2}-\delta)+o(1) > z(\frac{3}{2})-\xi +o(1)~.$$ For the
upper bound $S(\Gand[0]^{\tau_0})\leq z(\frac{3}{2})-o(1)$ we are
required to examine the second moment of
$S\left(\Gand[0](\frac{3}{2}+\delta)\right)$. Assume by
contradiction that:
\begin{equation}
  \label{eq-contradiction-s-tau0}
  \Pr\left[S\left(\Gand[0]^{\tau_0}\right) > z(\frac{3}{2}) + \xi \right] >
  \alpha \mbox{ for some fixed }\alpha,\xi > 0~,
\end{equation}
and choose $\delta>0$ small enough such that:
\begin{equation}\label{eq-delta-choice}
\left(\frac{z(\frac{3}{2})}{\xi/2}\right)^2
\left(\mathrm{e}^{12\delta}-\mathrm{e}^{2\delta}\right) < \alpha
~\mbox{ , and: }~
\left(\mathrm{e}^{4\delta}-1\right)z(\frac{3}{2})\leq \xi/2~.
\end{equation}
Set $T_0 = \left(\frac{3}{2}-\delta\right)\frac{n}{2}$ and
$T_1=\tau_0-1$, and note that, with high probability, $\Delta :=
\frac{T_1-T_0}{n/2}$ satisfies $\delta/2 \leq \Delta \leq 2\delta$,
and therefore we may assume that this holds. We consider
$\Gand[0]^T$ for $T=T_0,\ldots,T_1$, and let $S_0 =
S\left(\Gand[0]^{T_0}\right)$. As $I(\Gand[0]^T)
> 0$, the calculation which yielded \eqref{eq-z-K=0-eq} gives:
$$\mathbb{E} \left( S(\Gand[0]^{T+1}) ~|~S(\Gand[0]^T)\right) = S(\Gand[0]^T) + \frac{2}{n}\left(\frac{2
(S(\Gand[0]^T)-1)}{2-I(\Gand[0]^T)}+1\right)~, $$ and hence:
$$\left(1+\frac{2}{n}\right)S(\Gand[0]^T) \leq
\mathbb{E} \left( S(\Gand[0]^{T+1}) ~|~S(\Gand[0]^T)\right) \leq
\left(1+\frac{4}{n}\right)S(\Gand[0]^T) ~.$$ Therefore:
$$\mathrm{e}^{\Delta}S_0 \leq \mathbb{E}S(\Gand[0]^{\tau_0}) \leq \mathrm{e}^{2\Delta}
S_0~.$$ Similarly, we can write the expression for the second moment
of $S(\Gand[0]^{T_1})$:
$$\mathbb{E} \left( S(\Gand[0]^{T+1})^2 ~|~S(\Gand[0]^T)^2\right) \leq S(\Gand[0]^T)^2 \left(1+\frac{4}{n}
\cdot \frac{3}{2-I(\Gand[0]^T)}\right)~, $$ and hence:
$$\mathbb{E} \left( S(\Gand[0]^{\tau_0})^2\right) \leq
\mathrm{e}^{6\Delta} S_0^2~.$$ We obtain that:
$$\mathrm{Var}S(\Gand[0]^{\tau_0}) \leq \left(\mathrm{e}^{6\Delta}-\mathrm{e}^\Delta \right)S_0^2 \leq
\left(\mathrm{e}^{12\delta}-\mathrm{e}^{2\delta} \right)S_0^2~,$$
where the last inequality is by the fact that $\Delta \leq 2\delta$
almost surely. By \eqref{eq-delta-choice} and the fact that $S_0 <
z(\frac{3}{2})$ almost surely, this implies that: \begin{equation}
  \label{eq-var-s-tau0-bound}
  \mathrm{Var}S(\Gand[0]^{\tau_0}) < \alpha \left(\xi/2\right)^2~.
\end{equation}
By \eqref{eq-delta-choice} we have:
$$ \Pr\left[S\left(\Gand[0]^{\tau_0}\right) > z(\frac{3}{2}) + \xi
\right] \leq \Pr\left[S\left(\Gand[0]^{\tau_0}\right) -
\mathbb{E}\Gand[0]^{\tau_0} >
z(\frac{3}{2})\left(1-\mathrm{e}^{4\delta}\right) + \xi \right] \leq
$$
$$ \leq \Pr\left[|S\left(\Gand[0]^{\tau_0}\right) -
\mathbb{E}\Gand[0]^{\tau_0}| > \xi/2 \right]~,$$ and thus combining
Chebyshev's inequality with \eqref{eq-var-s-tau0-bound} gives:
$$ \Pr\left[S\left(\Gand[0]^{\tau_0}\right) > z(\frac{3}{2}) + \xi
\right] < \alpha~,$$ contradicting the assumption
\eqref{eq-contradiction-s-tau0}.

\begin{figure}
\centering
\includegraphics{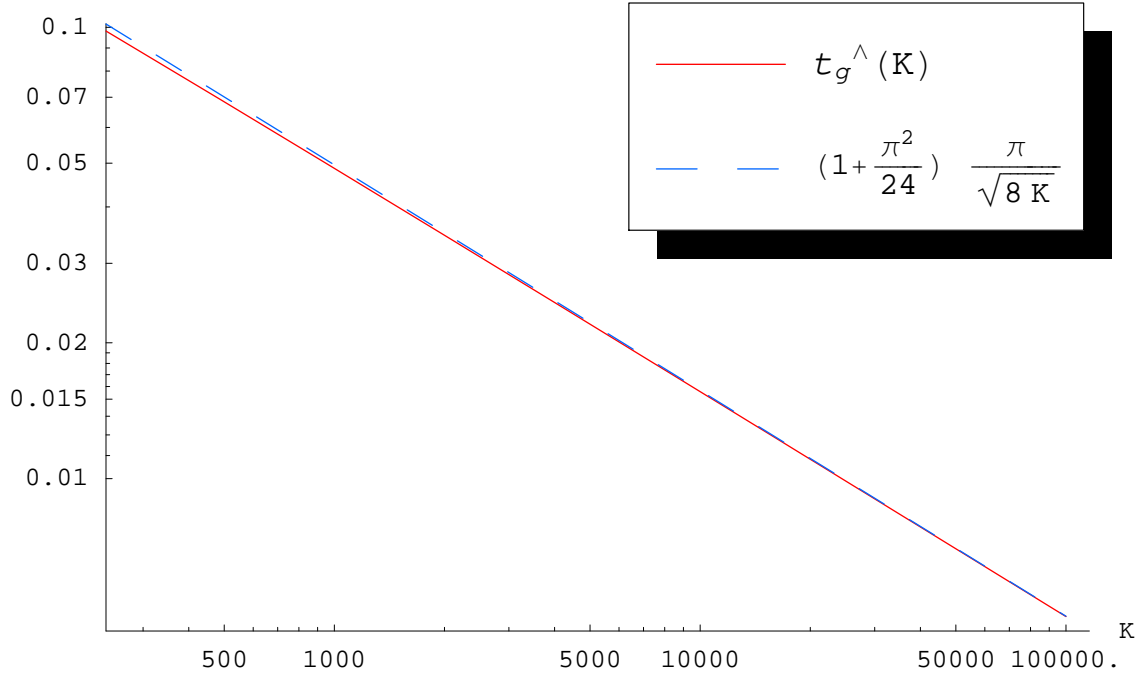}
\caption{Comparison of the values of $\tgand(K)$, obtained by
numerical solutions of the ODEs (\ref{y-eq})-(\ref{z-eq}), and the
asymptotic approximation of Theorem \ref{thm-K-0-K-gg-1}.
Logarithmic scale was used in both axes.} \label{fig::sing-approx}
\end{figure}

It remains to show that for $K \gg 1$, $\tgand(K) = (1+o(1))
\frac{\pi}{2\sqrt{2}}\left(1+\frac{\pi^2}{24}\right)\frac{1}{\sqrt{K}}$.
Consider equation \eqref{y-eq}; the initial condition $y(0)=1$
suggests that we examine the function $u(t)=1-y(t)$, which satisfies
$u(t)\ll 1$ near the origin:
$$-u' = \frac{u-1}{u^2(K-1)+1}~,$$ and substituting the initial condition of $y$ we have:
\begin{equation}\label{eq-u-ode}
u'=\frac{1-u}{u^2(K-1)+1}~,~u(0)=0~.
\end{equation}
As $u\ll 1$, the $u^2$-term at the denominator is negligible, and we
have: \begin{equation}\label{eq-u-ode-app} u' \approx \frac{1}{K
u^2+1}~,~u(0)=0~,
\end{equation}
where here and in what the follows, we denote by $\approx$ equality
up to leading order terms. Rearranging the last equation to the form
$Ku^2 u' + u' \approx 1$ and integrating, we obtain that, up to
leading order terms, $u$ satisfies the equation:
\begin{equation}\label{eq-u-cubic}
\frac{K}{3}u^3+u\approx t~.
\end{equation}
We note that this immediately gives $u\approx t$ for $u\ll
\frac{1}{\sqrt{K}}$, however we are interested in the behavior of
$u$ precisely when $t=\Theta(\frac{1}{\sqrt{K}})$. Applying
Caradano's solution to the above cubic equation gives the following
approximation of $u$ when $t\ll 1$ (and hence $u \ll 1)$:
$$u(t) \approx
\left(\frac{3t}{2K}+\sqrt{\frac{1}{K^3}+\left(\frac{3t}{2K}\right)^2}\right)^{1/3}+
\left(\frac{3t}{2K}-\sqrt{\frac{1}{K^3}+\left(\frac{3t}{2K}\right)^2}\right)^{1/3}
~,~ t\ll 1~.$$ Moving on to $z(t)$, we substitute $w=z-1$ in
equation \eqref{z-eq} and obtain the following:
$$w' = \frac{(w+1)^2+(K-1)(w+u)^2}{1+Ku^2 -u^2}\approx
\frac{(w+1)^2+(K-1)(w+u)^2}{1+K u^2}~.$$ Next, we may replace $w+1$
by $1$ whenever $w \ll 1$, and furthermore, whenever $w = \Omega(1)$
clearly the dominant term is $K(w+u)^2$. Altogether, we obtain the
following uniform approximation for $w$:
$$ w' \approx \frac{K(w+u)^2+1}{Ku^2+1} \approx \left(K(w+u)^2+1\right)u'~.$$
Adding $u'$ to both sides of the equation and rearranging, we
obtain:
$$\frac{w'+u'}{K(w+u)^2+2}\approx u'~,$$
hence if we define $v=\sqrt{\frac{K}{2}}(w+u)$ we obtain:
$$\frac{v'}{v^2+1}\approx \sqrt{2K}u'~,$$
and thus:
$$v \approx \tan(\sqrt{2K}u)~.$$
Returning to $z$, we get:
$$z \approx 1+ \sqrt{\frac{2}{K}}\tan(\sqrt{2K}u)-u~.$$
This implies that the singularity point $x_c$ satisfies
$\sqrt{2K}u(x_c) = \frac{\pi}{2}$. Recalling equation
\eqref{eq-u-cubic}, we have:
$$x_c \approx \frac{K}{3}u(x_c)^3 + u(x_c) =
\frac{K}{3}\frac{\pi^3}{16K\sqrt{2K}} + \frac{\pi}{2\sqrt{2K}} =
\frac{\pi}{2\sqrt{2}}\left(1+\frac{\pi^2}{24}\right)\frac{1}{\sqrt{K}}~.$$
Figure \ref{fig::sing-approx} shows the excellent agreement between
the above asymptotic approximation of $\tgand(K)$, and its value as
obtained by numerically solving the ODEs \eqref{y-eq} and
\eqref{z-eq} by Mathematica.  \qed

\noindent \textbf{Acknowledgement} The authors wish to thank Noga
Alon for helpful discussions and comments.

\end{document}